\documentclass{amsart}

\newtheorem{theorem}{Theorem}[section]

\theoremstyle{definition}
\newtheorem{definition}[theorem]{Definition}

\theoremstyle{remark}

\numberwithin{equation}{section}

\newcommand{\be}{\begin{equation}}
\newcommand{\ee}{\end{equation}}

\newcommand{\bea}{\begin{eqnarray}}
\newcommand{\eea}{\end{eqnarray}}
\newcommand{\Bea}{\begin{eqnarray*}}
\newcommand{\Eea}{\end{eqnarray*}}

\def\Cal{\mathcal}

\def\gnk{G_{n,k}}
\def\cgnk{\Cal G_{n,k}}

\def\ft{\hat f(\tau)}

\def\rn{\bbr^n}

\def\irn{\intl_{\bbr^n}}

\def\bbr{{\Bbb R}}

\def\bbh{{\Bbb H}}

\def\cos{{\hbox{\rm cos}}}

\def\gnk{G_{n,k}}

\def\part{\partial}
\def\intl{\int\limits}

\def\Gam{\Gamma}

\def\a{\alpha}

\def\vp{\varphi}

\def\sig{\sigma}
\def\lam{\lambda}
\def\z{\zeta}

\def\t{\tau}

\begin{document}

\title{On Y. Nievergelt's inversion formula for the Radon transform}

\author{E. OURNYCHEVA}

\address{University of Pittsburgh at Bradford,
 300 Campus Drive, 16701, Bradford, PA, USA }
  \email{elo10@pitt.edu}

\author{B. RUBIN}

\address{Department of Mathematics, Louisiana State
University \\ Baton Rouge, LA, 70803, USA}
\email{borisr@math.lsu.edu}

\thanks{ The second author was supported  by
 the  NSF grants
PFUND-137 (Louisiana Board of Regents) and DMS-0556157.}


\date{August 4 , 2009}



\begin{abstract}
We generalize Y. Nievergelt's inversion method for the Radon
transform on lines in the 2-plane to the $k$-plane  Radon transform
of continuous and $L^p$ functions on $\bbr^n$ for all $1\le k<n$.
\end{abstract}

\maketitle

\section{Introduction}

Inversion formulas for Radon transforms of different kinds are of
great importance in
  mathematics and its applications; see, e.g.,
\cite{Ehr, GGG, Hel, Mar, Na, Pal, QCK, RK,  Ru2, Str}, and
references therein. Since many of them are pretty involved,
especially for new-comers in the area, or applicable under essential
restrictions, every ``elementary" inversion method deserves special
consideration. In 1986 Yves Nievergelt came up with intriguing paper
\cite{Ni}, entitled ``Elementary inversion of Radon's transform".
His result
  can be stated as follows.

\begin{theorem} \label{nit} Let \be\label{ni2}
G_a(t)=\left \{
\begin{array} {ll} 1/(\pi a^2) & \mbox{ if   $\;|t|\le a$}, \\
{} \\
\displaystyle{\frac{1}{\pi a^2} \left
(1-\frac{1}{\sqrt{1-a^2/t^2}}\right )}   &\mbox{ if $\;|t|> a$;}
\end{array}
\right.
 \ee
$a>0$. Any continuous compactly supported function on the $2$-plane
can be reconstracted from the Radon transform over lines in this
plane by the formula \be\label{ni1} f(x,y)=\lim\limits_{a\to 0}\;
\frac{1}{\pi} \int_0^\pi\int_{-\infty}^\infty (Rf)(t-x\cos\,
\a-y\sin\a, \a)\, G_a(t)\, dt\, d\a,\ee where the double integral on
the right-hand side equals the average of $f$ over the disc of
radius $a$ centered at $(x,y)$.
\end{theorem}

Formulas (\ref{ni2}) and (\ref{ni1}) indeed look elementary. The
following questions arise:

1. What is the basic idea of the Nievergelt's method from the point
of view of modern developments?

2. Is this method  applicable in the same  elementary form to
$k$-plane Radon transforms on $\bbr^n$ for all $1\le k<n$ and
arbitrary continuous or $L^p$ functions, satisfying minimal
assumptions at infinity?

In the present article we answer these questions and indicate
possible generalizations.

\vskip 0.3truecm

\noindent{\bf Notation and main results.} Let $\cgnk$ and $\gnk$ be
the {\it affine} Grassmann manifold of all non-oriented $k$-planes
$\t$ in $\rn$, and the ordinary Grassmann manifold
 of  $k$-dimensional  subspaces $\z$ of $\rn$, respectively.  Given $\zeta
\in \gnk$, each vector $x \in \rn$ can be written as $x = (x^\prime,
x'') = x^\prime+x''$ where $x^\prime\in\zeta$ and
$x''\in\zeta^\perp$,
  $\zeta^\perp$ being the orthogonal complement to $\zeta $ in $\rn$.
 Each $k$-plane $\tau$  is parameterized by the pair
$( \zeta, x'')$ where $\zeta \in \gnk$ and $ x'' \in \zeta^\perp$.
 The manifold  $\cgnk$ will be endowed with the product measure $d\t=d\zeta dx''$,
where $d\z$ is the
 $SO(n)$-invariant measure  on $\gnk$ of  total mass
$1$, and $dx''$ denotes the usual volume element on $\zeta^\perp$.
We  write $C_0=C_0 (\bbr^n)$ for the space of continuous functions
on $\bbr^n$ vanishing at infinity; $\sigma_{n-1}  = 2 \pi^{n/2}/
 \Gamma(n/2)$ denotes the area of the unit sphere $S^{n-1}$ in $\rn$.

The $k$-plane transform   of a function $f$  on $\rn$ is a function
$\hat f$ on  $\cgnk$ defined by  \be \ft  = \intl_\zeta f(x^\prime +
x'') \,dx^\prime, \qquad \t=(\zeta, x'') \in \cgnk. \ee This
expression is finite  for all $\t\in \cgnk$ if $f$ is continuous and
decays like $O(|x|^{-\lam})$ with $ \lam>k$. Moreover \cite{Ru2, So,
Str},
 if $f\in L^p (\rn)$, $1\le p<n/k$, then  $\ft$ is finite for almost
all planes $\t\in\cgnk$. The above-mentioned bounds for $\lam$ and
$p$ are best possible.

Following \cite{Ru1}, we define the wavelet-like transform \be\label
{wav} (W^*_a
 \vp)(x)=\frac{1}{a^n}\intl_{\cgnk}  \vp (\t) \,w \Big (\frac {|x
 -\t|}{a}\Big ) \, d\t,  \qquad a>0,\ee
where $|x
 -\t|$ denotes the Euclidean distance  between the point $x\in \rn$ and the $k$-plane
$\t$.
\begin{theorem} \label {tcon1} \cite [Th. 3.1]{Ru1} Let $ f \in L^p, \; 1 \le
p <n/k$, and let $\psi (|\cdot|)$ be a radial function on $\rn$,
which has  an integrable  decreasing radial majorant. If $w$ is a
solution of the Abel type integral equation \be\label{(3.3)}  c
r^{2-n} \! \intl_0^r \! s^{n-k-1} w(s) (r^2 \! - \! s^2)^{k/2 -1}\,
ds\! = \! \psi (r), \quad c \! = \! \frac{\sig_{k-1} \,
\sig_{n-k-1}}{\sig_{n-1}},
 \ee
 then \be \label {34} (W^*_a \hat f)(x)=
\irn f(x-ay)\,\psi (|y|)\, dy, \ee
 and therefore, \be\label{inv}
\lim\limits_{a \to 0}\,(W^*_a \hat f)(x)= \lam  f(x),\qquad
\lam=\irn \psi (|x|) \,dx.\ee The limit in (\ref{inv}) is understood
in the $L^p$-norm and in the almost everywhere sense.  If $f \in C_0
\cap L^p$ for some $1 \le p <n/k$, then (\ref{inv}) holds uniformly
on $ \rn$.
\end{theorem}

 This theorem is a core of the convolution-backprojection method for the $k$-plane
Radon transform, and the most difficult task is to choose relatively
simple functions $w$ and $\psi$ satisfying (\ref{(3.3)}). The crux
 is that the left-hand side of (\ref{(3.3)}) has, in general,
  a bad behavior when $r\!\to\! \infty$. Hence, {\it to achieve integrability of $\psi$, the solution $w$ must be
  sign-changing.}

Our first observation is that the essence of Y. Nievergelt's Theorem
\ref {nit} can be presented in the language of Theorem
 \ref{tcon1} as follows.
\begin{theorem} \label{nit1} Let $k=1$, $ n=2$,
\be\label{psi0} \psi (r)= \left \{
\begin{array} {ll} 1 & \mbox{ if   $\;0\le r\le 1$}, \\
0  &\mbox{ if $\;r> 1$.}
\end{array}
\right.
 \ee
Then (\ref{(3.3)}) has a  solution \be\label{ni2w} w(r)=\left \{
\begin{array} {ll} 1 & \mbox{ if   $\;0\le r\le 1$}, \\
{} \\
\displaystyle{1-\frac{r}{\sqrt{r^2-1}}}
&\mbox{ if $\;r> 1$,}
\end{array}
\right.
 \ee
  such that for every compactly supported continuous function $f$
on  $\bbr^2$,
 \be\label{(niev)} (W^*_a \hat
f)(x)=\intl_{|y|<1} f(x-ay) \, dy, \ee where $W^*_a$ is
 the  wavelet-like
 transform (\ref{wav}) generated by $w$.
\end{theorem}

For the convenience of  presentation, we will keep to the following
 convention.
\begin{definition}
 The convolution-backprojection algorithm  in Theorem
\ref{tcon1} will be called
 {\it the Nievergelt's method } if $\psi$  is chosen according to (\ref{psi0}).
\end{definition}

Of course, Theorem  \ref{tcon1} deals with essentially more general
classes of functions than Theorem \ref{nit}, however, the main
 focus of our article is different: {\it we want to find auxiliary functions $\psi$ and $w$, having
possibly simple analytic expression.}

\vskip 0.3truecm

\noindent{\bf Theorem A.} \hfil

\noindent {\rm (i)}  {\it The Nievergelt's method is applicable to
the
 $X$-ray transform (the case $k=1$) in any dimension. Namely, if
$\psi$ is chosen according to (\ref{psi0}), then
 (\ref{(3.3)}) has a solution  \be\label{w-k1a} w(r)=\left \{
\begin{array} {ll} \!1 & \!\mbox{ if   $\;0\le r\le 1$}, \\
{} \\

  \!\displaystyle{-\frac{  \Gam\Big(\frac{n-1}{2}\Big )\, r^{3-n} }{2 \sqrt{\pi}
  \,\Gam\Big(\frac{n}{2}\Big)} \!\intl_0^{1}  v^{n/2-1} \,(r^2\!-\!v)^{-3/2}\,dv} \! &\mbox{ if
$\;r>1$,}
\end{array}
\right.
 \ee
and  inversion formula (\ref{inv}) holds.

 \noindent {\rm (ii)} If $k>1$, then Nievergelt's method is
inapplicable. }

An integral in (\ref{w-k1a}) can be expressed through the
hypergeometric function and explicitly evaluated in some particular
cases; see Section 2.2. For instance, if $n=2$, then (\ref{w-k1a})
is  the  Nievergelt's function (\ref {ni2w}).

To include all $1\le k<n$, we modify
 the Nievergelt's method  by
choosing $\psi (r)$ in a different way as follows.

\vskip 0.3truecm

\noindent {\bf Theorem B.} {\it Let \be\label{2ellq} \psi(r)=\left
\{
\begin{array} {ll} 0 & \mbox{ if   $\;0\le r\le 1$}, \\
{} \\
\displaystyle{\frac{(r^2-1)^\ell}{r^{n+2\ell +1}} } &\mbox{ if $\;r>
1$;}
\end{array}
\right.
 \ee
$\ell \ge 0$. Then the corresponding function $\psi(|\cdot|)$ has a
decreasing radial majorant  in $L^1(\rn)$ and
 (\ref{(3.3)}) has the following  solution:

\vskip 0.3truecm

\noindent {\rm (i)} In the case $k=2\ell; \; \ell=1,2, \ldots \, :$
\be\label{2ellz} w(r)=\left \{
\begin{array} {ll} 0 & \mbox{ if   $\;0\le r\le 1$}, \\
{} \\
\displaystyle{c_\ell\,r^{2+2\ell -n}\,\left (\frac{1}{2r}\,
\frac{d}{dr}\right )^\ell \left [\frac{(r^2-1)^\ell}{r^{2\ell
+3}}\right ]} &\mbox{ if $\;r> 1$,}
\end{array}
\right.
 \ee
$$
c_\ell=\frac{\Gam (n/2-\ell)}{\Gam (n/2)}. $$

 \vskip 0.3truecm

\noindent {\rm (ii)} In the case $k=2\ell+1; \; \ell=0,1,2, \ldots
\, :$ \be\label{2ellc} w(r)\!=\!\left \{
\begin{array} {ll}\! 0 & \mbox{ if   $\;0\!\le \!r\!\le \!1$}, \\
{} \\
\displaystyle{\!\tilde c_\ell\,r^{3+2\ell -n}\,\left (\frac{1}{2r}\,
\frac{d}{dr}\right )^{\ell+1} \left
[\frac{(r^2\!-\!1)^{\ell+1/2}}{r^{2\ell +2}}\right ]} &\mbox{ if
$\;r> 1$,}
\end{array}
\right.
 \ee
 $$
\tilde c_\ell=\frac{\Gam ((n-1)/2-\ell)\, \ell !}{\Gam (n/2)\, \Gam
(\ell +3/2)}. $$ In both cases the inversion result in Theorem \ref
{tcon1} is valid.}

\vskip 0.3truecm

Theorems {\bf A} and  {\bf B}  are proved in Sections 2 and 3,
respectively.

\vskip 0.3truecm

\noindent{\bf Possible generalizations.} \hfill

\vskip 0.3truecm

\noindent $1^\circ$. The convolution-backprojection method is
well-developed in the general context of totally geodesic Radon
transforms on spaces of constant curvature. Apart of $\bbr^n$, the
latter include the $n$-dimensional unit sphere $S^n$  and the
hyperbolic space $\bbh^n$; see \cite {BR, Ru9, Ru1}. As above, the
key role in this theory belongs to a certain Abel type integral
equation and the relevant sign-changing solution $w$. Moreover,
passage to the limit in (\ref{inv}) as $a \to 0$,
 can be replaced by integration in  $a$ from $0$
to $\infty$ against the dilation-invariant measure $da/a$. This
leads to
 inversion formulas, which resemble the classical Calder\'on's
 identity
 for continuous wavelet transforms \cite{FJW}. The corresponding wavelet
 function is determined as a solution of a similar Abel type integral
equation; see \cite {BR, Ru9, Ru1} for details. In all these cases
analogues of Theorems {\bf A} and {\bf B}
 can be obtained. We leave this exercise to the interested reader.

\vskip 0.3truecm

\noindent $2^\circ$.  Unlike the classical $k$-plane transforms on
$\bbr^n$, the corresponding transforms on matrix spaces \cite {GK,
OOR, OR1, OR2} are much less investigated. To the best of our
knowledge, no pointwise
 inversion formulas (i.e., those, that do not contain operations in the sense
  of distributions) are available for these transforms if
 the latter are applied  to arbitrary continuous or $L^p$ functions  $(p\neq
 2)$. One of the reasons of our interest in
Nievergelt's idea
 is that it might be applicable to the matrix case. Moreover, as in $1^\circ$,
 it may pave the way to implementation of wavelet-like transforms in the corresponding
 reconstruction formulas.  We plan to
 study these questions in our forthcoming publication.

\section{The case $k=1$}

\setcounter{equation}{0}

\subsection{Proof of Theorem A}\label{2.1x}

We will be dealing with Riemann-Liouville fractional integrals
\be\label{ril} (I^\a_{a+} g)(u) = \frac{1}{\Gamma(\a)}\int^u_a
(u-v)^{\a-1} g (v) \,dv, \qquad \a>0.\ee
 Changing variables, we transform the basic integral equation
\be \label{basic} c r^{2-n} \! \intl_0^r \! s^{n-k-1} w(s) (r^2 \! - \! s^2)^{k/2
-1}\, ds\! = \! \psi (r), \quad c \! = \! \frac{\sig_{k-1} \,
\sig_{n-k-1}}{\sig_{n-1}},
 \ee
(cf. (\ref{(3.3)})) to the form \be I^{k/2}_{0+} \tilde w =\tilde
\psi,\ee where \be \tilde w(u)=u^{(n-k)/2-1} w(\sqrt u),\qquad
\tilde \psi(u)=\frac{2u^{n/2-1}}{c\, \Gam (k/2)}\, \psi(\sqrt u).\ee
Suppose that $\psi$ is defined by  (\ref{psi0}). If $0<u\le 1$, then,
by homogeneity, \be\label{xxx} \tilde w(u)=c'u^{(n-k)/2-1}, \qquad
c'=\frac{2\Gam (n/2)}{c\, \Gam (k/2)\,\Gam ((n\!-\!k)/2)}=1, \ee $c$
being the constant from (\ref{basic}).
 Hence, for $u>1$, we necessarily have
\be\label{uuu} (I^{k/2}_{1+} \tilde w)(u) =-\frac{1}{\Gam
(k/2)}\int_0^1 (u-v)^{k/2 -1}v^{(n-k)/2-1}\, dv. \ee If $k\geq 2$,
this equation has no solution $\tilde w\in L^1_{loc}(1,\infty)$,
because, otherwise, we get $$\lim\limits_{x\to 1^+}(l.h.s)=0, \qquad
\lim\limits_{x\to 1^+}(r.h.s)=const\neq 0.
$$
This proves the second statement in Theorem {\bf A}.

Consider the case $k=1$. If  $0<u\le 1$, then $ \tilde
w(u)=u^{(n-1)/2-1}$. If $u>1$, then, setting $w_0(u)=u^{(n-1)/2-1}$,
from  (\ref{uuu}) we get \bea (I^{1/2}_{1+} [\tilde w
-w_0])(u)&=&-\frac{1}{\sqrt \pi}\int_0^u (u-v)^{-1/2}
v^{(n-1)/2-1}\, dv\nonumber\\&=&-c_n\, u^{n/2-1}, \qquad
c_n=\frac{\Gam ((n-1)/2)}{\Gam (n/2)}.\nonumber\eea This gives
$$
\tilde w(u)=w_0(u)-\frac{c_n}{\sqrt \pi}\,\frac{d}{du}\int_1^u
(u-v)^{-1/2} v^{n/2-1}\, dv$$ or (set $\int_1^u=\int_0^u -
\int_0^1$)
$$
\tilde w(u)=-\frac{c_n}{2\sqrt \pi}\,\int_0^1(u-v)^{-3/2}
v^{n/2-1}\, dv.$$ Thus, if  $I^{1/2}_{0+} \tilde w=\tilde \psi$,
then, necessarily, \be \label {ppp}\tilde w(u)= \left \{
\begin{array} {ll} u^{(n-1)/2-1} & \mbox{ if   $0<u\leq 1$}, \\
{} \\
  \displaystyle{-\frac{c_n}{2\sqrt \pi}\,\int_0^1(u-v)^{-3/2} v^{n/2-1}\,
dv}  &\mbox{ if $u>1$.}
\end{array}
\right.
 \ee
One can readily see that function (\ref{ppp}) is locally integrable
on $(0,\infty)$. Let us prove that it  satisfies $(I^{1/2}_{0+}
\tilde w)(u)=\tilde \psi (u)$ for all $u>0$. It suffices to show
 that $(I^{1/2}_{0+} \tilde w)(u)\equiv 0$ when $u>1$.  We have $I^{1/2}_{0+} \tilde
 w=I_1 -I_2$, where
\bea
 I_1&=&\frac{1}{\sqrt \pi}\int_0^1(u-v)^{-1/2} v^{(n-1)/2-1}\,
 dv,\nonumber\\
I_2&=&\frac{c_n}{2\pi}\int_1^u (u-v)^{-1/2} \,dv\int_0^1(v-s)^{-3/2}
s^{n/2-1}\,
 ds.\nonumber\eea
Both integrals can be expressed in terms of hypergeometric
functions. For $I_1$, owing to   3.197 (3) and 9.131 (1) from
\cite{GRy}, we obtain \bea
 I_1&=&\frac{u^{-1/2}}{\sqrt \pi}\, B\left(\frac{n-1}{2},1\right
 )\, F\left(\frac{1}{2},\frac{n-1}{2}; \frac{n+1}{2}; \frac{1}{u}\right
 )\nonumber\\
&=&\label {888}\frac{\Gam ((n-1)/2)}{\sqrt {\pi}\,\Gam
((n+1)/2)}\,(u-1)^{-1/2}\, F\left(1,\frac{1}{2}; \frac{n+1}{2};
\frac{1}{1-u}\right
 ).\eea
For  $I_2$, changing the order of integration and using  \cite[3.238
(3)]{GRy}, we have
$$
I_2=\frac{c_n\,(u-1)^{1/2}}{\pi}\int_0^1\frac{
s^{n/2-1}\,(1-s)^{-1/2}}{u-s}\, ds.$$
 By
\cite[3.228 (3)]{GRy} this expression coincides with  (\ref{888}).

To complete the proof, we recall that $w(r)=r^{3-n}\tilde w(r^2)$,
which gives
$$
 w(r)= \left \{
\begin{array} {ll} 1 & \mbox{ if   $0\le r\leq 1$}, \\
{} \\
  \displaystyle{-\frac{c_n\,r^{3-n}}{2\sqrt \pi}\,\int_0^1(r^2-v)^{-3/2} v^{n/2-1}\,
dv}  &\mbox{ if $r>1$.}
\end{array}
\right.
 $$
This coincides with (\ref{w-k1a}).

\subsection{Examples}

Let us give some examples of functions $w(r)$  defined by
(\ref{w-k1a}) in the case $r>1$.  By \cite [3.197(3)]{GRy},
$$
w(r)=-\frac{\Gam ((n-1)/2)\, r^{-n}}{2\sqrt {\pi}\, \Gam (n/2+1)}\,
F\left (\frac{3}{2}, \,\frac{n}{2}; \,\frac{n}{2}+1;\,
\frac{1}{r^2}\right ).$$ Keeping in mind that $F(a,b; c; z)=F(b,a;
c; z)$, and using formulas 156, 203, and 211 from
\cite[7.3.2]{PBM3}, we obtain:

\vskip 0.3truecm

 \noindent For $n=2$:   \be\label{w-k1-n2}
w(r)=-\frac{1}{2r^2}\,F\left (\frac{3}{2},\, 1; \,2;\,
\frac{1}{r^2}\right )=1-\frac{r}{\sqrt{r^2-1}}. \nonumber\ee

 \noindent For $n=3$:  \be\label{w-k1-n3} w(r)\!=\!-\frac{ 1}{2r^3\sqrt
{\pi}\, \Gam (5/2)}\,F\left (\frac{3}{2},
\,\frac{3}{2};\,\frac{5}{2}; \,\frac{1}{r^2}\right
)\!=\!\frac{2}{\pi} \Big(\arcsin
\frac{1}{r}\!-\!\frac{1}{\sqrt{r^2\!-\!1}} \Big). \nonumber\ee

 \noindent For $n=4$:   \be\label{w-k1-n4} w(r)=-\frac{
1}{8r^4}\,F\left (\frac{3}{2},\, 2; \,3;\,
 \frac{1}{r^2}\right )=
1-\frac{2r^2-1}{2r\sqrt{r^2-1}}. \nonumber \ee

\section{The general case}

\setcounter{equation}{0}

As in Section \ref {2.1x}, our main concern is  integral equation
(\ref{basic}), which is equivalent to \be\label{vvv} I^{k/2}_{0+}
\tilde w =\tilde \psi,\ee  \be\label{jjj} \tilde w(u)=u^{(n-k)/2-1}
w(\sqrt u),\qquad \tilde \psi(u)=\frac{2u^{n/2-1}}{c\, \Gam (k/2)}\,
\psi(\sqrt u).\ee We want to find relatively simple functions $w$
and $\psi$, which are admissible in the basic Theorem \ref  {tcon1}
and such that the corresponding functions $\tilde w$ and
$\tilde\psi$ obey (\ref{vvv}). It is convenient to consider the
cases of $k$ even and $k$ odd separately.

\subsection{The  case of $k$ even } Let $k=2\ell; \; \ell=1,2, \ldots \, $. We choose
\be\label{2ell5} \psi(r)=\left \{
\begin{array} {ll} 0 & \mbox{ if   $\;0\le r\le 1$}, \\
{} \\
\displaystyle{\frac{(r^2-1)^\ell}{r^{n+2\ell +1}} }
&\mbox{ if $\;r> 1$.}
\end{array}
\right.
 \ee
The corresponding function $x\to \psi(|x|)$ obviously has  a
decreasing radial majorant  in $L^1(\rn)$. By  (\ref{jjj}), equation
(\ref{vvv}) is equivalent to \be\label{eee} (I^{\ell}_{1+} \tilde
w)(u) =\frac{c_\ell\, (u-1)^\ell}{u^{\ell +3/2}}, \qquad
c_\ell=\frac{\Gam (n/2-\ell)}{\Gam (n/2)}, \quad u>1.\ee The
$\ell$th derivative \be\label{nnn} \tilde w(u) =\left
[\frac{c_\ell\, (u-1)^\ell}{u^{\ell +3/2}}\right ]^{(\ell)}=
c_\ell\, \sum\limits_{j=0}^\ell  {\ell \choose j}\,
[(u-1)^\ell]^{(j)}\, [u^{-\ell-3/2}]^{(\ell-j)}\ee is integrable on
$(1,\infty)$. The fact that  (\ref{nnn}) satisfies (\ref{eee}) can
be easily checked using integration by parts. Thus, the pair of
functions $ w(r)$ and $ \psi(r)$, defined by \be\label{2ell}
w(r)=\left \{
\begin{array} {ll} 0 & \mbox{ if   $\;0\le r\le 1$}, \\
{} \\
\displaystyle{c_\ell\,r^{2+2\ell -n}\,\left (\frac{1}{2r}\, \frac{d}{dr}\right )^\ell
\left [\frac{(r^2-1)^\ell}{r^{2\ell +3}}\right ]}
&\mbox{ if $\;r> 1$,}
\end{array}
\right.
 \ee
and (\ref{2ell5}), falls into the scope of Theorem \ref  {tcon1} and
the ``even part" of Theorem {\bf B} is proved.

\subsection{The  case of $k$ odd } Let $k=2\ell+1; \; \ell=0,1,2, \ldots \, $.
We  define $\psi(r)$ by (\ref{2ell5}), as above. Then, instead of
(\ref{eee}), we have \be\label{eee1} (I^{\ell+1/2}_{1+} \tilde w)(u)
=\frac{c'_\ell\, (u-1)^\ell}{u^{\ell +3/2}}, \qquad
c'_\ell=\frac{\Gam ((n-1)/2-\ell)}{\Gam (n/2)}, \quad u>1.\ee This
gives $\tilde w(u) = g^{(\ell+1)}(u)$, where
$$
g(u)=\frac{c'_\ell}{\sqrt \pi}\intl_1^u
\frac{(u-v)^{-1/2}(v-1)^\ell}{v^{\ell +3/2}}\, dv=\tilde
c_\ell\,\frac{(u-1)^{\ell+1/2}}{u^{\ell +1}},
$$
$$
\tilde c_\ell=\frac{\Gam ((n-1)/2-\ell)\, \ell !}{\Gam (n/2)\, \Gam
(\ell +3/2)},$$ (use  \cite[3.238(3)]{GRy}). Let us show that
$\tilde w(u) \equiv g^{(\ell+1)}(u)$ satisfies (\ref{eee1}).
Integrating by parts, we have $I^{\ell+1/2}_{1+} \tilde w
=I^{1/2}_{1+} g'$, where
$$
g'(u)=\frac{dg(u)}{du}=\tilde c_\ell\, \left [\frac{(\ell +1/2)\,
(u-1)^{\ell-1/2}}{u^{\ell +1}}- \frac{(\ell +1)\,
(u-1)^{\ell+1/2}}{u^{\ell +2}}\right ].
$$
This gives $I^{\ell+1/2}_{1+} \tilde w =\tilde c_\ell\,(I_1-I_2)$
where (use  \cite[3.238(3)]{GRy} again) \bea
I_1&=&\frac{\ell+1/2}{\sqrt \pi}\intl_1^u \frac{(u-v)^{-1/2}(v-1)^{\ell-1/2}}{v^{\ell +1}}\, dv=\frac{\Gam (\ell +3/2)}{\ell !}\, \frac{(u-1)^{\ell}}{u^{\ell+1/2}}, \nonumber\\
I_2&=&\frac{\ell+1}{\sqrt \pi}\intl_1^u \frac{(u-v)^{-1/2}(v-1)^{\ell+1/2}}{v^{\ell +2}}\, dv=\frac{\Gam (\ell +3/2)}{\ell !}\, \frac{(u-1)^{\ell+1}}{u^{\ell+3/2}}.\nonumber\eea
Hence,
$$
I_1-I_2=\frac{\Gam (\ell +3/2)}{\ell !}\, \frac{(u-1)^{\ell}}{u^{\ell+3/2}}
$$
and therefore,
$$
(I^{\ell+1/2}_{1+} \tilde w)(u) =\frac{\tilde c_\ell\,\Gam (\ell
+3/2)}{\ell !}\, \frac{(u-1)^{\ell}}{u^{\ell+3/2}}=c'_\ell\,
\frac{(u-1)^{\ell}}{u^{\ell+3/2}},
$$
as desired. Thus,  functions $ w(r)$ and $ \psi(r)$, defined by
\be\label{2ell} w(r)=\left \{
\begin{array} {ll} 0 & \mbox{ if   $\;0\le r\le 1$}, \\
{} \\
\displaystyle{\tilde c_\ell\,r^{3+2\ell -n}\,\left (\frac{1}{2r}\,
\frac{d}{dr}\right )^{\ell+1} \left
[\frac{(r^2-1)^{\ell+1/2}}{r^{2\ell +2}}\right ]} &\mbox{ if $\;r>
1$,}
\end{array}
\right.
 \ee
and (\ref{2ell5}), obey  Theorem \ref  {tcon1}. This completes the
proof of Theorem {\bf B}.

\bibliographystyle{amsplain}

\end{document}